\documentclass[11pt]{article}

\usepackage{times}
\usepackage{amsthm}
\usepackage{amsmath}
\usepackage{amssymb}
\usepackage{mathrsfs}

\newcommand\ie{{\em i.e.}~}

\def\A{{\mathcal A}}

\def\th{\theta}
\def\Th{\Theta}

\def\N{\mathbb{N}}

\def\X{\mathscr X}
\def\BC{\mathsf{BC}}

\def\QQ{\mathbf Q}

\def\C{\mathcal{C}}

\def\B{{\mathcal B}}

\def\H{\mathcal H}

\def\Q{\mathcal Q}
\def\P{\mathcal P}
\def\R{\mathcal{R}}

\def\BB{\mathfrak B}
\def\CC{{\mathfrak C}}
\def\PP{{\mathfrak P}}
\def\RR{{\mathfrak R}}

\def\s{{\mathsf s}}
\def\J{\mathcal J}
\def\JJ{\mathfrak J}

\def\I{{\rm 1\kern-.26em I}}

\def\si{\sigma}
\def\Si{\Sigma}

\def\sp{\mathop{\mathrm{sp}}\nolimits}
\def\spe{{\rm sp}_{{\rm ess}}}

\def\d{\mathrm{d}}

\def\1{\mathfrak{1}}
\def\0{\mathfrak{0}}

 \def\hb{\hbar}

\def\<{\langle}
\def\>{\rangle}

\def\Op{\mathfrak{Op}}

\providecommand{\CC}{\mathfrak{C}}

\def\supp{\mathop{\mathrm{supp}}\nolimits}
\def\ltwo{\mathsf{L}^{\:\!\!2}}

\newtheorem{Theorem}{Theorem}[section]
\newtheorem{Remark}[Theorem]{Remark}
\newtheorem{Lemma}[Theorem]{Lemma}

\newtheorem{Corollary}[Theorem]{Corollary}
\newtheorem{Proposition}[Theorem]{Proposition}
\newtheorem{Definition}[Theorem]{Definition}
\newtheorem{Example}[Theorem]{Example}

\numberwithin{equation}{section}

\begin{document}


\title{Rieffel's pseudodifferential calculus and\\ spectral analysis of quantum Hamiltonians}

\date{\today}

\author{M. M\u antoiu 
}
\date{\small}
\maketitle \vspace{-1cm}


\begin{abstract}
We use the functorial properties of Rieffel's pseudodifferential calculus to study families of operators associated to topological dynamical systems acted by a symplectic space. Information about the spectra and the essential spectra are extracted from the quasi-orbit structure of the dynamical system. The semi-classical behavior of the families of spectra is also studied.
\end{abstract}

\section{Introduction}\label{duci}

In \cite{Rie1}, Marc Rieffel significantly extended the core the Weyl pseudodifferential calculus. His main purpose was to provide
a unified framework for a large class of examples in deformation quantization (cf. also \cite{Rie2}). The emerging formalism has very nice functorial properties, which makes it virtually an efficient tool in other
directions. In this article we intend to apply it to certain problems in spectral analysis.

A way to summarize (a restricted version of) \cite{Rie1} is to say that it gives an exact contravariant functor from the category of locally compact dynamical systems with group $\Xi:=\mathbb R^{2n}$ to a category of (usually non-commutative) $C^*$-algebras, also endowed with an action of $\Xi$. To achieve this, the canonical symplectic form on $\mathbb R^{2n}$ plays an important role; it contributes to defining a composition law deforming the point-wise multiplication of functions acting on the space of the dynamical system. The resulting non-commutative $C^*$-algebras are essentially composed of functions that can be interpreted as observables of some quantum systems. Under suitable circumstances, the real-valued ones can be represented as bounded self-adjoint operators in Hilbert spaces and they might have physical meaning. The main theme of the present article is to show that many spectral properties of these operators can be tracked back to properties of the underlying topological dynamical system, just by using the way Rieffel quantization has been constructed. In particular, the quasi-orbit structure of the dynamical system as well as the nature of each quasi-orbit will play central roles.

One of our starting points was a question asked by Vladimir Georgescu in connection with the problem of determining the
essential spectrum of anisotropic differential and pseudodifferential operators. We are going to outline it in a
framework which is convenient subsequently. But other applications of Rieffel's formalism to spectral analysis are
contained, as explained later on.

We set $\X:=\mathbb R^n$ and $\Xi:=\X\times\X^*\cong\mathbb R^{2n}$, where $\X^*$ is the dual of the space $\X$, the
duality being denoted simply by $\X\times\X^*\ni(x,\xi)\mapsto x\cdot\xi$. To suitable functions $h$ defined on "the
phase space" $\Xi$, one assigns by "quantization" operators acting on function $u:\X\to\mathbb C$ by
\begin{equation}\label{op}
\left[\mathfrak{Op}(h)u\right](x):=(2\pi)^{-n}\int_\X\int_{\X^*}dx\,d\xi\,e^{i(x-y)\cdot
\xi}\,h\left(\frac{x+y}{2},\xi\right)u(y).
\end{equation}
This is basically the Weyl quantization and, under convenient assumptions on $h$, \eqref{op} makes sense and has nice
properties in the Hilbert space $\H:=\ltwo(\X)$ or in the Schwartz space $\mathcal S(\X)$.

Let $h:\Xi\to \mathbb R$ be an elliptic symbol of strictly positive order $m$. This means that $h$ is smooth
and satisfies estimates of the form
\begin{equation}\label{sm}
\left\vert\left(\partial^\alpha_x\partial_\xi^\beta h\right)(x,\xi)\right\vert\le C_{\alpha,\beta}(1+|\xi|)^{m-|\beta|},\ \ \ \ \ \forall
\alpha,\beta\in\N^n,\ \ \ \forall (x,\xi)\in\Xi
\end{equation}
and
\begin{equation}\label{el}
|h(x,\xi)|\ge C(1+|\xi|)^m,\ \ \ \ \ \forall (x,\xi)\in\Xi,\ \ |\xi|\ {\rm large\ enough}.
\end{equation}
It is well-known that under these assumptions $\Op(h)$ makes sense as an unbounded self-adjoint operator in $\H$,
defined on the $m$'th order Sobolev space. The problem is to evaluate the essential spectrum of this operator; it
comes out that the relevant information is contained in the behavior at infinity of $h$ in the $x$ variable.

This one is conveniently taken into account through an Abelian algebra $\A$ composed of uniformly continuous functions
un $\X$, which is invariant under translations (if $\varphi\in\A$ and $y\in \X$ then
$\theta_y(\varphi):=\varphi(\cdot+y)\in\A$). Let us also assume (for simplicity) that $\A$ is unital and contains the
ideal $\C(\X)$ of all complex continuous functions on $\X$ which converge to zero at infinity. We ask that the
elliptic symbol $h$ of strictly positive order $m$ also satisfy
\begin{equation}\label{circu}
\left(\partial^\alpha_x\partial^\beta_\xi h\right)(\cdot,\xi)\in\A,\ \ \ \ \ \forall\,\alpha,\beta\in\mathbb N^n,\ \forall\,\xi\in\X^*.
\end{equation}
Then the function $h$ extends continuously on $\Omega\times\X^*$, where $\Omega$ is the Gelfand spectrum of the
$C^*$-algebra $\A$ and it is a compactification of the locally compact space $\X$. By translational invariance, it is a
compact dynamical system under an action of the group $\X$. After removing the orbit $\X$, one gets a $\X$-dynamical
system $\Omega_\infty:=\Omega\setminus\X$; its quasi-orbits (closure of orbits) contain the relevant information about the
essential spectrum of the operator $H:=\Op(h)$. For each quasi-orbit $\Q$, one constructs a self adjoint operator
$H_\Q$. It is actually the Weyl quantization of the restriction of $h$ to $\Q\times\X^*$, suitably reinterpreted.
Using the notations $\sp(T)$ and $\spe(T)$, respectively, for the spectrum and the essential spectrum of an operator
$T$, one gets finally
\begin{equation}\label{sp}
\spe(H)=\overline{\bigcup_\Q\sp(H_\Q)}.
\end{equation}
Many results of this kind exist in the literature, some of them for special type of functions $h$, but with less
regularity required, others including anisotropic magnetic fields and others formulated in a more geometrical framework. We only cite \cite{ABG,Da,GI1,GI2,GI3,HM,LMN,LN,LMR,LS,Ma1,MPR2,RRS}; see also references therein. As V. Georgescu remarked (\cite{GI1,GI2} and private communication), when the
function $h$ does not diverge for $\xi\to\infty$, the approach is more difficult and should also take into
account the asymptotic values taken by $f$ in "directions contained in $\X^*$". One of our goals is to achieve this in
a sufficiently general framework.

A very efficient tool for obtaining some of the results cited above was the crossed product, associated to
$C^*$-dynamical systems. In the setting presented before, one uses the action $\theta$ of $\X$ by translations on the
$C^*$-algebra $\A$ to construct a larger, non-commutative $C^*$-algebra $\A\rtimes_\theta\X$. After a partial Fourier
transform, this one can be seen to be generated by pseudodifferential operators of strictly negative order, with
coefficients in $\A$. So it will contain resolvent families of elliptic strictly positive order Weyl operators
satisfying \eqref{circu} and the structure of the crossed product will rather easily imply spectral results. A basic
fact is that the crossed product is a functor, also acting on equivariant morphisms, and that it behaves nicely with
respect to quotients and direct sums. One drawback is, however, that $\xi$-anisotropy cannot be treated easily. The
symbols of order $0$ are not efficiently connected to the crossed products (treating them as multiplier would not be enough
for our purposes).

To overcome this, we are going to use the general pseudodifferential calculus of \cite{Rie1}. It is strong tool,
containing as a particular case the crossed product construction connected with strictly negative order Weyl operators. It
has as basic data the action $\Theta$ of a vector space (as our "phase space" $\Xi$) on a $C^*$-algebra $\B$ (even a
non-commutative one), together with a skew-symmetric  linear operator $J:\Xi\to\Xi$ that serves to twist the
product on $\B$. This twisting is done first on the set of smooth elements of $\B$ under the action. Then a $C^*$-norm
is found on the resulting non-commutative $^*$-algebra. The outcome will be a new $C^*$-algebra $\mathfrak B$ (the
quantization of $\B$, composed of pseudodifferential symbols) also endowed with an action of the vector space $\Xi$.
In \cite{Rie1} it is shown that one gets a strict deformation quantization of a natural Poisson structure defined on
$\B$ by the couple $(\Theta,J)$, but this will not concern us here. It will be more important to note the functorial
properties of the construction; they will be reviewed in Section \ref{sectra}, in a particular setting suited to our purposes. In Section \ref{secintra} this mechanism is used to define ideals and quotients associated to quasi-orbits, as a basis for the forthcoming proofs concerning spectral analysis.

In order to get a $C^*$-norm on the set of smooth vectors with the new, non-commutative product, Rieffel uses Hilbert
module techniques. He deliberately neglects exploring Schr\"odinger-type representations. In the most general setting
they might not be available, and when they are they could be unfaithful, which is a serious drawback for his aims.
For us, Schr\"odinger representations are essential: they make the necessary connection to the operators
we would like to study, and they also offer tools of investigation. So we dedicate a section to their
definition under the assumption that the initial $C^*$-algebra is Abelian, so it defines canonically a locally compact topological space $\Si$ (the Gelfand spectrum). We make use of the quasi-orbit structure of the dynamical system $(\Sigma,\Theta,\Xi)$ to arrive in more
familiar spaces, directly connected to $\Xi$, where the traditional pseudodifferential theory works. For a fixed classical observable $f$, the outcome is a family of operators $\{H_\si\}_{\si\in\Si}$ acting in the Hilbert space $\H:=\ltwo(\mathbb R^n)$, indexed by the points of the space $\Si$ and grouped together in classes of unitary equivalence along the orbits under the action $\Th$.
Each individual operator might be more complicated than a usual pseudodifferential operator in $\mathbb R^n$; this is connected to the fact that the action $\Th$ is a general one. From this point of view, the set up is interesting even if the orbit structure of $(\Si,\Th,\Xi)$ is poor, as in the case of topologically transitive systems for which one of the orbits is dense.

In Section \ref{etrata} we start our spectral analysis for the operators $H_\si$, using the formalism
presented before. We get first a spectral inclusion results connected to the
hyerarquisation of representations by the quasi-orbit structure. It follows that two operators $H_{\si_1}$ and $H_{\si_2}$ will be equi-spectral if the points $\si_1$ and $\si_2$ generates the same quasi-orbit. This is weaker in general that the property of belonging to the same orbit, which would imply that  $H_{\si_1}$ and $H_{\si_2}$ are even unitarily equivalent.
We also include a decomposition formula, used subsequently in the analysis of essential spectra.

In Section \ref{etratata} we present our results on the essential spectrum of pseudodifferential operators with $(x,\xi)$
(phase-space) anisotropy and defined by general actions of phase-space. The proofs exploit rather straightforwardly the properties of Rieffel's quantization, a
simple strategy to convert structural information about $C^*$-algebras into spectral results on operators
naturally connected to them and some lemmas about dynamical systems proved in Section \ref{secintra}. It is shown that the operators associated to a
certain type of points, called of the second kind, have no discrete spectrum.  This happens in particular for those belonging to a minimal quasi-orbit. In the opposite case (points $\si$ of the first kind), the situation is more interesting. The essential spectrum of $H_\si$ is the closed union of spectra of operators $H_{\si'}$ associated to the non-generic points $\si'$ (those belonging to the quasi-orbit generated by $\si$, but generating themselves strictly smaller quasi-orbits).

A section is dedicated to some examples, illustrating mainly Theorem \ref{mailatata}.

Then we turn to a random setting,
defined by an invariant ergodic probability on $\Si$. In Section \ref{secintre} information about the a.e. constancy of the spectrum is once again deduced from the formal properties of Rieffel's calculus
and from arguments in ergodic theory. It is also shown that with probability one the operators $H_\si$ have void discrete spectrum. Such results (and many others) are quite standard in the theory of random Hamiltonians. We included this short section because we can give rather precise statements and the proofs uses nicely the previous framework.

Rieffel's calculus also contains a deformation parameter $\hb$, which in some situations can be assimilated to
Planck's constant. In the limit $\hb\to 0$ one recovers the initial classical data (codified in the form of a Poisson algebra) from the deformed
structures. In the present article, almost everywhere, the value $\hb=1$ is fixed. In Section \ref{serata}, reintroducing $\hb$ in the formalism, one gets among others families $\left(H^\hb\right)_{\hb\in[0,1]}$ of
Hamiltonians defined by a symbol $f$.  Then we show that the family $\left\{\sp\left(H^\hb\right)\right\}_{\hb\in(0,1]}$
of spectra of the quantum observables converges for $\hb\to 0$ to the spectrum $\sp\left(H^0\right)=\overline
{f(\Si)}$ of the classical observable.

We stress that most of our spectral results do not use the functional calculus, so they stand for general elements of the relevant $C^*$-algebra and not only for self-adjoint ones.

Up to our knowledge, most of the results contained in this paper are new, at least in this form and for this class of Hamiltonians. But any expert in spectral analysis for quantum Hamiltonians would like to see the analog of these results for unbounded symbols and, maybe, for a version including magnetic fields. To achieve this, the technique of affiliation of unbounded observables to $C^*$-algebras (\cite{ABG,GI1,GI2}) can be used. But in order to give affiliation a wide applicability and a deep theoretical foundation, Rieffel's calculus should also be extended in two directions: First, it should include unbounded elements, connected to but not contained in $C^*$ algebras. Second, it should incorporate group $2$-cocycles much more complicated than the one defined by the canonical symplectic form on $\Xi$ (thus generalizing the magnetic pseudo-differential calculus developed and applied to spectral problems in \cite{MP1,MPR2,IMP1,LMR}). This is work for the future.

\section{Pseudodifferential operators \`a la Rieffel}\label{setra}

\subsection{Rieffel's pseudodifferential calculus}\label{sectra}

We shall recall briefly some constructions and results from \cite{Rie1}. Whenever our aims allow it, we'll choose to
simplify; most noteworthy, the initial (un-quantized) algebra will be {\it Abelian} and the vector space will be endowed with a {\it non-degenerate} bilinear anti-symmetric form. Some convention will also be different. For the moment we fix Planck's constant setting $\hb=1$, but we shall come back to this point in Section \ref{serata}.

The initial object, containing {\it the classical data}, is a quadruplet
$\left(\Sigma,\Theta,\Xi,[\![\cdot,\cdot,]\!]\right)$. $\left(\Xi,[\![\cdot,\cdot]\!]\right)$ is a $d$-dimensional
symplectic vector space. The number $d$ is pair and there is no loose of generality to imagine that a Lagrangean decomposition $\,\Xi=\X\times\X^*$ was given, with $\X^*$ the dual of the $n$-dimensional vector space $\X$, and that for
$X:=(x,\xi),\,Y:=(y,\eta)\in\Xi$, the symplectic form reads
\begin{equation}\label{sympl}
[\![X,Y]\!]:=x\cdot\eta-y\cdot \xi.
\end{equation}

A continuous action $\Theta$ of $\Xi$ by homeomorphisms of the locally compact space $\Sigma$ is also given. For
$(\sigma,X)\in\Sigma\times\Xi\,$ we are going to use all the notations
$$
\Theta(\sigma,X)=\Theta_X(\sigma)=\Theta_\sigma(X)\in\Sigma
$$
for the $X$-transformed of the point $\sigma$. The function $\Theta$ is continuous and the homeomorphisms $\Th_X,\Th_Y$ satisfy $\Th_X\circ\Th_Y=\Th_{X+Y}$ for every $X,Y\in\Xi$.

We denote by $\C(\Sigma)$ the Abelian $C^*$-algebra of all complex continuous functions on $\Sigma$
that are arbitrarily small outside large compact subsets of $\Sigma$. When $\Sigma$ is compact,
$\C(\Si)$ is unital. The action $\Theta$ of $\Xi$ on $\Sigma$ induces an action of $\Xi$ on
$\C(\Si)$ (also denoted by $\Theta$) given by
\begin{equation}\label{indu}
\Theta_X(f):=f\circ\Theta_X.
\end{equation}
This action is strongly continuous, \ie for any $f\in\C(\Sigma)$ the mapping
\begin{equation}\label{ciuca}
\Xi\ni X\mapsto\Theta_X(f)\in\C(\Sigma)
\end{equation}
is continuous. We denote by $\C^\infty(\Sigma)$ the set of elements $f\in\C(\Sigma)$ such that the
mapping \eqref{ciuca} is $C^\infty$; it is a dense $^*$-algebra of $\C(\Sigma)$.
It is also a Fr\'echet algebra for the family of semi-norms
$$
|f|_k:=\sum_{|\alpha|\le k}
\big\Vert\partial_X^\alpha\big(\Theta_X(f)\big)\big|_{X=0}\big\Vert_{\C(\Sigma)}\,,\qquad k\in\N.
$$

To quantize the above structure, one introduces on $\C^\infty(\Sigma)$ the product
\begin{equation}\label{rodact}
f\,\#\,g:=\pi^{-2n}\int_\Xi\int_\Xi dYdZ\,e^{2i[\![Y,Z]\!]}\,\Theta_Y(f)\,\Theta_Z(g),
\end{equation}
suitably defined by oscillatory integral techniques and set simply $f^*(\si):=\overline{f(\si)},\,\forall\si\in\Si$.
One gets a $^*$-algebra $\big(\C^\infty(\Sigma),\#,^*\big)$, which
admits a $C^*$-completion $\CC(\Si)$ in a $C^*$-norm $\Vert\cdot\Vert_{\CC(\Si)}$.

The action $\Theta$ leaves $\C^\infty(\Sigma)$ invariant and extends to a strongly continuous action on
the non-commutative $C^*$-algebra $\CC(\Sigma)$; the space $\CC^\infty(\Si)$ of $C^\infty$-vectors in
$\CC(\Sigma)$ coincides with $\C^\infty(\Sigma)$.

Actually the quantization transfers to $\Xi$-morphisms, and this will be crucial in the sequel. Let
$\left(\Sigma_j,\Theta_j,\Xi,[\![\cdot,\cdot]\!]\right)$, $j=1,2$, be two classical data with Abelian $C^*$-algebras $\C(\Sigma_j)$ and let $\R:\C(\Si_1)\to\C(\Si_2)$ a $\Xi$-morphism, \ie a ($C^*$-)morphism intertwining the
two actions $\Theta_1,\Theta_2$. Then $\R$ acts coherently on $C^\infty$-vectors and extends to a morphism
$\RR:\CC(\Si_1)\to\CC(\Si_2)$ that also intertwines the corresponding actions. In this way, one obtains a
covariant functor between two categories of $C^*$-algebras endowed with $\Xi$-actions, the
algebras being Abelian in the first category.

The functor is exact: it preserves short exact sequences of $\Xi$-morphisms. Namely, if $\J$ is a (closed,
self-adjoint, two-sided) ideal in $\C(\Si)$ that is invariant under $\Theta$, then its quantization $\JJ$ can be
identified with an invariant ideal in $\CC(\Si)$ and the quotient $\CC(\Si)/\JJ$ is canonically isomorphic to the
quantization  of the quotient $\C(\Si)/\J$  under the natural quotient action. Composing Rieffel's functor with the Gelfand functor, we get a contravariant functor from the category of locally compact $\Xi$-dynamical systems to a category of non-commutative $C^*$-dynamical systems with group $\Xi$.

An important example is given by {\it$\Xi$-algebras}, \ie $C^*$-algebras $\B$ composed of
bounded, uniformly continuous function on $\Xi$, under the additional assumption that the action $\mathcal T$ of $\Xi$ on
itself by translations, raised to functions, leaves $\B$ invariant. Consequently, by Gelfand theory, there exists a
continuous function $:\Xi\mapsto\Sigma$ with dense image, which is equivariant with respect to the actions
$\mathcal T$ on $\Xi$, respectively $\Theta$ on $\Sigma$. The function is injective if and only if $\C(\Xi)\subset\B$.

The largest such $C^*$-algebra $\B$ is $\BC_{{\rm u}}(\Xi)$, consisting of all the
bounded uniformly continuous functions $:\Xi\mapsto\mathbb C$. It coincides with the family of functions
$
g\in\BC(\Xi)$ (just bounded and continuous) such that
$$
\Xi\ni X\mapsto g\circ\mathcal T_X=g(\cdot+X)\in\BC(\Xi)
$$
is continuous. Then the Fr\'echet $^*$-algebra of $C^\infty$-vectors is
$$
\BC_{{\rm u}}(\Xi)^\infty\equiv\BC^\infty(\Xi):=\{f\in C^\infty(\Xi)\mid |\left(\partial^\alpha
f\right)(X)|\le C_\alpha,\,\forall\,\alpha,X\}.
$$
It might be illuminating to note that it coincides with $S^0_{0,0}(\Xi)$, one of H\"ormander's symbol classes.

Another important particular case is $\B=\C(\Xi)$ (just put $\Sigma=\Xi$ in the general construction). It is shown in \cite{Rie1} that at the quantized level one gets the usual Weyl calculus and the emerging non-commutative $C^*$-algebra $\CC(\Xi)$ is
isomorphic to the ideal of all compact operators on an infinite-dimensional separable Hilbert space.

\subsection{$C^*$-algebras associated to quasi-orbits}\label{secintra}

It is convenient to have a closer look at the quasi-orbit structure of the dynamical system $(\Si,\Theta,\Xi)$ in connection with $C^*$-algebras and representations.

For each $\si\in\Si$, we write $E_\si:=\overline{\Theta_\si(\Xi)}$ for {\it the quasi-orbit generated by} $\si$ and set
$$
\P_\si:\C(\Sigma)\to\BC_{\rm u}(\Xi),\quad\P_\sigma(f):=f\circ\Theta_\sigma\,.
$$
The range of the $\Xi$-morphism $\P_\si$ is called $\B_\sigma$ and it is a $\Xi$-algebra.
Defining analogously $\P'_\si:\C(E_\sigma)\to\BC_{\rm u}(\Xi)$ one gets a $\Xi$-monomorphism with the
same range $\B_\sigma$, which shows that the Gelfand spectrum of $\B_\si$ can be identified with the quasi-orbit $E_\si$.

For each quasi-orbit $E$, one has the natural restriction map
$$
\R_E:\C(\Si)\to\C(E),\quad\R_E(f):=f|_E,
$$
which is a $\Xi$-epimorphism. Actually one has $\P_\si=\P'_\si\circ \R_{E_\si}$.

Being respectively invariant under the actions $\Theta$ and $\mathcal T$, the $C^*$-algebras $\C(E)$
and $\B_\sigma$ are also subject to Rieffel deformation. By quantization, one gets $C^*$-algebras and
morphisms
$$
\RR_E:\CC(\Si)\to\CC(E),\qquad\PP_\si:\CC(\Si)\to\BB_\si,\qquad\PP'_\si:\CC(E_\si)\to\BB_\si,
$$
satisfying $\PP_\si=\PP'_\si\circ \RR_{E_\si}$. While $\RR_E$ and $\PP_\si$ are epimorphisms, $\PP'_\si$ is an isomorphism.

We denote by $\QQ(\Si)\equiv\QQ(\Si,\Theta,\Xi)$ the family of all the
quasi-orbits. For every $E\in\QQ(\Si)$, the restriction of the action (also denoted by $\Theta$) defines a dynamical
subsystem $(E,\Theta,\Xi)$. If $E$ is a quasi-orbit in $(\Si,\Theta,\Xi)$ we set
$$
\QQ(E):=\{F\in\QQ(\Si)\mid F\subset E\}\qquad\hbox{and}\qquad\QQ_0(E):=\QQ(E)\setminus\{E\}.
$$
For any $E\in\mathbf Q(\Si)$ let us denote by $\C^E(\Si)$ the $C^*$-subalgebra of $\C(\Sigma)$
composed of elements that vanish on the closed invariant set $E$. Obviously, it is an invariant ideal
coinciding with the kernel of the morphism $\R_E$, and the quotient $\C(\Sigma)/\C^E(\Sigma)$ can
be identified with $\C(E)$. By applying Rieffel's functor one gets
$$
\CC(\Sigma)/\CC^E(\Sigma)\equiv\CC(\Sigma)/\ker(\RR_E)\cong\CC(E).
$$
More generally, for any two closed subsets $\Sigma_1\subset\Sigma_2$ of $\Sigma$, we denote by
$\C^{\Sigma_1}(\Sigma_2)$ the closed ideal of $\C(\Sigma_2)$ composed of functions which are zero
on $\Sigma_1$.

Some points of a quasi-orbit generate the quasi-orbit, others do not. We write $\si\propto E$ if $E=E_\si$ (this may be stronger than $\si\in E$). Every quasi-orbit $E$ decomposes as an invariant disjoint union $E^{\rm g}\sqcup E^{\rm n}$, where
$E^{\rm g}:=\{\si\propto E\}$ is the dense set of {\it generic points} and the complement
$E^{\rm n}=\bigcup_{F\in\mathbf Q_0(E)}F$ is composed of {\it non-generic points}.

When $E^{{\rm n}}=\varnothing$, \ie when all the orbits contained in $E$ are dense, one says that $E$ is {\it minimal};
the points generating a minimal compact quasi-orbit are called {\it almost-periodic}. Clearly $E$ is minimal iff it does not contain non-trivial invariant closed subsets and also iff $\C(E)$ is $\Xi$-{\it simple},
\ie it does not contain non-trivial invariant ideals.

\medskip
An important role in our subsequent spectral analysis will be played by the following notions:

\begin{Definition}\label{ghel}
\begin{enumerate}
\item
A $C^*$-subalgebra $\B$ of $\,\BC(\Xi)$ is {\rm of
the first kind} if $\,\C(\Xi)\subset\B$ and it is {\rm of the second kind} if $\,\B\cap\C(\Xi)=\{0\}$.
\item
This can be applied to $\Xi$-algebras $\B$ and, in particular, to those of the form $\B_\sigma=\P_\si[\C(\Si)]$ for
some $\si\in\Si$. Accordingly, $\sigma\in \Sigma$ (and also the quasi-orbit $E_\si$) is {\rm of the first} (resp. {\rm second) kind} if $\,\B_\sigma$
is of the first (resp. second) kind.
\end{enumerate}
\end{Definition}

In general $\,\B$ might fail to be of one of the two kinds, but
this is not possible for $\,\Xi$-algebras. We are grateful to Serge Richard for this remark.

\begin{Lemma}\label{serge}
Any $\,\Xi$-algebra $\B$ is either of the first or of the second kind.
\end{Lemma}

\begin{proof}
We know from \cite[Thm.~5.1]{dLM60} that $\C(\Xi)$ does not have non-trivial
translational-invariant $C^*$-subalgebras (this is based essentially on the Stone-Weierstrass Theorem). Thus $\B\cap\C(\Xi)=\{0\}$ or $\B\cap\C(\Xi)=\C(\Xi)$, which implies the claim.
\end{proof}

Consequently, we get the invariant decomposition $\Si=\Si_I\sqcup\Si_{I\!I}$, where the points in $\Si_I$ are of the first kind and those of $\Si_{I\!I}$ of the second kind.

We will also use later the next result:

\begin{Lemma}\label{nimalu}
Let $\B$ be a $\Xi$-algebra with Gelfand spectrum $E$. If $E$ is minimal then either $\B=\C(\Xi)$ or $\B$ is of the second kind.
\end{Lemma}

\begin{proof}
Obviously $E=\Xi$ is minimal, so we can exclude this case.

If $\B$ is not of the second kind, then it is of the first kind by Lemma \ref{serge}.
But then $\C(\Xi)$ is a non-trivial invariant ideal of $\B$, which contradicts the
minimality of $E$.
\end{proof}

\begin{Lemma}\label{ciudat}
If there exists an open orbit $\mathcal O$ in the quasi-orbit $E$, then this one is dense, it is the single dense orbit and it coincides with the set of all generic elements of $E$.
\end{Lemma}

\begin{proof}
To see this, let $\mathcal O'$ be another orbit. Then $\overline{\mathcal O'}\cap\mathcal O=\varnothing$, so $\mathcal O'$
cannot be dense. Since a dense orbit exists, this one must be $\mathcal O$. Thus all the elements of $\mathcal O$ are generic, and those not belonging to $\mathcal O$ are not, since they belong to some $\mathcal O'\ne\mathcal O$.
\end{proof}

The next result will be useful subsequently in the study of essential spectra:

\begin{Proposition}\label{detai}
For $\si\in\Si$ of the first kind, let us set $\mathcal O\equiv\mathcal O_\si=\Th_\si(\Xi),\,E\equiv E_\si:=\overline{\mathcal O}$ and $\B:=\B_\si$ and consider the isomorphism
$\P'_\si:\C(E)\to \B\subset\BC_{{\rm u}}(\Xi)$. Then $\mathcal O$ is open and coincides with $E^{{\rm g}}$ and
\begin{equation}\label{detoi}
(\P'_\si)^{-1}\left[\C(\Xi)\right]=\C^{E^{{\rm n}}}(E).
\end{equation}
\end{Proposition}

\begin{proof}
If $\B$ is unital we leave things as they are, if not we embed $\B$ into its minimal unitization $\widehat\B$ and identify the Gelfand spectrum $\,\widehat B$ of $\,\widehat\B\,$ with the Alexandrov compactification $B\sqcup\{\infty\}$ of the Gelfand spectrum $B$ of $\B$. We are going to treat the non-unital case; the unital one needs less arguments.

Since the Gelfand spectrum
of $\BC(\Xi)$ is the Stone-\v Cech compactification $\beta\Xi$ of $\Xi$, the monomorphism
$\widehat\B\mapsto\C(\beta\Xi)$ induces a continuous surjection $\beta\s:\beta\Xi\to \widehat B$ which restricts
to a continuous map $\,\s:\Xi\to \widehat B$ with dense range. It is easy to see that $\s$ is injective, since we assumed that $\widehat \B$ contains $\C(\Xi)$; thus $\widehat B$ is {\it a compactification} of $\Xi$.

Clearly $\,\beta\s\left(\beta\Xi\setminus\Xi\right)\,$ is a closed subset of $\widehat B$, containing the point $\infty$. It follows that $\Xi$ can be identified to a dense {\it open} subset of the locally compact space $B$; the functions in $\B$ are characterized by the fact that they extend continuously from $\Xi$ to $B$.

We use now the fact that $\B$ is a $\,\Xi$-algebra. The action $\mathcal T$ of $\,\Xi$ on itself by translations extends to a topologically transitive dynamical system $(B,\mathcal T,\Xi)$ and the mapping $\Th_\si:\Xi\rightarrow\mathcal O\subset E$ extends to an isomorphism $\widetilde\Th_\si$ between the dynamical systems $(B,\mathcal T,\Xi)$ and $(E,\Th,\Xi)$, which defines by pull-back the $C^*$-isomorphism $\P'_\si$. With this picture in mind, it is clear that $\mathcal O$ is open; thus, by Lemma \ref{ciudat}, it coincides with $E^{{\rm g}}$.

Let $g\in\C(B)\equiv \B$. Then $g\in\C(\Xi)$ iff $g$ vanishes on $B\setminus \Xi$, which happens iff $\,\left(\P'_\si\right)^{-1}(g)=g\circ\widetilde\Th_\si^{-1}$ vanishes on $\widetilde\Th_\si(B\setminus\Xi)=E\setminus E^{{\rm g}}=E^{{\rm n}}$.
\end{proof}

\subsection{Representations and families of Hamiltonians}\label{secintraa}

We now construct representations of $\CC(\Si)$ using the points $\si$ of its Gelfand spectrum. On
$\BC^\infty(\Xi)$ one can apply the usual Schr\"odinger representation in $\H:=\ltwo(\X)$
\begin{equation}\label{opa}
\Op:\BC^\infty(\Xi)\to\mathbf B(\H)
\end{equation}
given by \eqref{op}, rigorously defined as an oscillatory integral. By the Calder\'on-Vaillancourt Theorem \cite{Fol89}, \eqref{opa} is a well-defined continuous function; it is also a $^*$-morphism with respect to complex conjugation and symbol composition.

For the sake of formalism, let us denote by $\C(\nu\Xi)$ the $\Xi$-algebra $\BC_{{\rm u}}(\Xi)$ (this is an awkward way to give a notation to its Gelfand spectrum - sometimes called {\it the uniform compactification of $\,\Xi=\mathbb R^{2n}$}) and by $\CC(\nu\Xi)$ its Rieffel quantization. The common set of smooth vectors is $\BC^\infty(\Xi)$, on which $\Op$ is a $^*$-morphism with respect to the structure of $\CC(\nu\Xi)$.

\begin{Proposition}\label{vai}
The mapping $\Op$ extends to a faithful representation of $\CC(\nu\Xi)$ in $\H$.
\end{Proposition}

\begin{proof}
For $\B=\C(\Xi)$ it is known that the quantization $\CC(\Xi)$ is isomorphic (essentially by a partial Fourier transform) to the crossed product $\C(\X)\rtimes_\tau\X$, with $\tau$ the action of $\X$ on $\C(\X)$ by translations (we recall that $\Xi=\X\times\X^*$). One can infer this from \cite{Rie1}, Example 10.5. Then it follows that $\Op$ realizes an isomorphism from $\CC(\Xi)$ to the ideal $\mathbf K(\H)$ of all compact operators in $\H$.

Let us recall that an ideal (always supposed closed and bi-sided) $\mathcal K$ in a $C^*$-algebra $\mathcal A$ is called {\it essential} if, for any $a\in\mathcal A$, from $a\mathcal K=\{0\}$ we deduce that $a=0$. Another equivalent condition is to have
$$
\parallel b\parallel_\mathcal A=\sup\{\parallel kb\parallel_\mathcal A\,\mid\, k\in\mathcal K,\,\parallel k\parallel_\mathcal A\equiv\parallel k\parallel_\mathcal K=1\},\ \ \ \ \forall\,b\in\mathcal A.
$$
Proposition 5.9 in \cite{Rie1} asserts that quantifying essential ideals one gets essential ideals. Now $\C(\Xi)$ is an essential ideal in $\C(\nu\Xi)$, so $\CC(\Xi)$ will be an essential ideal in $\CC(\nu\Xi)$. On the other hand, it is well-known and easy to prove that $\mathbf K(\H)$ is an essential ideal in $\mathbf B(\H)$. Thus, for $h\in\BC^\infty(\Xi)$,
we can write
$$
\parallel \Op(h)\parallel_{\mathbf B(\H)}=\sup\left\{\parallel K\Op(h)\parallel_{\mathbf B(\H)}\,\mid\,K\in\mathbf K(\H),\, \parallel K\parallel_{\mathbf B(\H)}=1\right\}=
$$
$$
=\sup\left\{\parallel \Op(k)\Op(h)\parallel_{\mathbf B(\H)}\,\mid\,k\in\CC(\Xi),\, \parallel \Op(k)\parallel_{\mathbf B(\H)}=1\right\}=
$$
$$
=\sup\left\{\parallel k\#h\parallel_{\CC(\Xi)}\,\mid\,k\in\CC(\Xi),\, \parallel k\parallel_{\CC(\Xi)}=1\right\}=\parallel h\parallel_{\CC(\nu\Xi)}.
$$
This is enough to prove the statement. We avoided to use the explicit definition of the $C^*$-norm on $\CC(\nu\Xi)$.
\end{proof}

For each $\Xi$-algebra $\B$, we restrict $\Op$ from $\BC^\infty(\Xi)$ to $\,\B^\infty=\BB^\infty$ (the dense
$^*$-algebra of smooth vectors of $\B\,$) and then we extend it to a faithful representation in $\H$ of
the $C^*$-algebra $\BB$.
We can apply the construction to the $C^*$-algebras $\BB_\si$. By composing, we get a family $\big\{\Op_\sigma:=\Op\circ\PP_\sigma\big\}_{\sigma\in\Sigma}$ of
representations of $\CC(\Sigma)$ in $\H$, indexed by the points of $\Sigma$. For $f\in\CC^\infty(\Sigma)$ one has $\PP_\sigma(f)\in\BB_\sigma^\infty=\B_\si^\infty$, and the action on $\H$ is given by
\begin{equation}\label{teisn}
\big[\Op_\sigma(f)u\big](x)=(2\pi)^{-n}\int_\X\d y\int_{\X^*}\d\xi\,
e^{i(x-y)\cdot\xi}f\Big[\Theta_{\left(\frac{x+y}2,\xi\right)}(\sigma)\Big]u(y)
\end{equation}
in the sense of oscillatory integrals. If the function $f$ is real, all the operators $\Op_\si(f)$ will be self-adjoint.

\begin{Remark}\label{vigonia}
{\rm The point $\si$ was called of the first kind when $\C(\Xi)\subset\B_\si$. In such a case $\mathbf K(\H)\subset\Op(\BB_\si)$, thus $\Op_\si$ is irreducible. Notice that $\Op_\si$ is faithful exactly when $\PP_\si$ is injective, i.e. when $\P_\si$ is injective, which is obviously equivalent to $E_\si=\Si$. Consequently, if the dynamical system is not topologically transitive, none of the Schr\"odinger-type representations $\Op_\si$ will be faithful. On the other hand we have the easy to prove identity}
$$
\parallel f\parallel_{\CC(\Si)}=\underset{\si\in\Si}{\sup}\parallel\Op_\si(f)\parallel_{\mathbf B(\H)},\ \ \ \ \ f\in\CC(\Si).
$$
\end{Remark}

\begin{Remark}\label{inbun}
{\rm Let us denote by $\mathcal S(\X)$ the Schwartz space on $\X=\mathbb R^n$ and by $\mathcal S(\Xi)$ the Schwartz space on $\Xi=\mathbb R^{2n}$. The corresponding duals are, respectively, the spaces of tempered distributions $\mathcal S^*(\X)$ and $\mathcal S^*(\Xi)$. But $\,\mathbf B(\H)$ is contained in $\mathbf B[\mathcal S(\X),\mathcal S^*(\X)]$ (the space of all linear and continuous operators $:\mathcal S(\X)\rightarrow\mathcal S^*(\X)$), and \cite{Fol89} the later is isomorphic to $\mathcal S^*(\Xi)$, by an extension of the representation $\Op$. Therefore, for every $\si\in\Si$, the $C^*$-closure $\,\BB_\si$ of $\BB^\infty=\B^\infty=\{f\circ\Th_\si\mid f\in\C^\infty(\Si)\}$ can be realized as a $C^*$-algebra of temperate distributions on $\Xi$. This gives a more concrete flavor to the spectral results of the next Section.}
\end{Remark}

To conclude part of the discussion above, a single element $\,f\in\CC(\Sigma)\,$ leads to a family
$\big\{H_\sigma:=\Op_\sigma(f)\big\}_{\sigma\in\Si}$ of bounded operators in $\H$ (in Quantum Mechanics we are mainly interested in the self-adjoint case). Such
families usually appear in disguise, as we are going to explain now.

Let $F:\Sigma\times\Xi\to\mathbb C$ be a function satisfying $F(\;\!\cdot\;\!,X)\in\C(\Sigma)$ for
each $X\in\Xi$, and such that $X\mapsto F(\;\!\cdot\;\!,X)\in\C(\Sigma)$ is smooth. In addition,
assume that $F$ satisfies {\it the equivariance condition}
$$
F\big(\Theta_Y(\sigma),X\big)=F(\sigma,X+Y)\quad\hbox{for all }\si\in\Si\,\hbox{ and }\,X,Y\in\Xi,
$$
which is very often imposed on physical reasons. Put $F_\sigma(X):=F(\si,X)$; then one interprets
$\big\{F_\sigma\big\}_{\sigma\in\Sigma}$ as an equivariant family of classical observables defined on the phase-space $\Xi$, that can
be transformed into quantum observables $\widetilde{H_\sigma}:=\Op(F_\sigma)$ by the usual Weyl
calculus. So, apparently, every operator $\widetilde{H_\sigma}$ has its own symbol $F_\sigma$.
Define $f:\Sigma\to\mathbb C$ by $f(\sigma):=F(\sigma,0)$. Then $f$ belongs to
$\,\C^\infty(\Sigma)\subset\CC(\Sigma)$. Moreover, one has
$$
[\P_\sigma(f)](X)
=(f\circ\Theta_\sigma)(X)
=F\big(\Theta_X(\sigma),0\big)
=F(\sigma,X)
=F_\sigma(X).
$$
Thus
$$
H_\sigma\equiv\Op_\sigma(f)=\Op(F_\sigma)\equiv\widetilde{H_\sigma},
$$
and we are in the framework presented above.

In other situations, a single operator $H$ is given as the Weyl quantization of a real function $f$ defined in phase space.
The behavior of $f$ requires the introduction of a $\Xi$-algebra $\B$ with Gelfand spectrum $\Si$; then $H=\Op(f)$ is
a represented version of $f$ seen as an element of the deformed algebra $\BB$. In favorable circumstances
$\Si$ is a compactification of the phase-space $\Xi$ on which $\Xi$ acts by homeomorphisms, and one has $H=\Op_\si(f)$ with $\si=0\in\Xi\subset\Si$. Other operators
$H_X:=\Op_X(f)$ defined by the points $X$ of the orbit $\Xi$ are all unitarily equivalent with $H$ (see Theorem \ref{colin}  for a more general statement). But the remaining
family $\big\{H_\si:=\Op_\si(f)\big\}_{\si\in\Si\setminus\Xi}$ is also useful in the spectral analysis of $H$. For instance, they give decompositions of the essential spectrum of the operator $H$; this is a particular case of Theorem \ref{mailatata}.

These remarks justify studying spectra of the family of operators $\{H_\sigma\}_{\sigma\in\Sigma}$ by using Rieffel's quantization.

\section{Spectral analysis for pseudodifferential operators}\label{sintra}
\setcounter{equation}{0}

\subsection{Spectra}\label{etrata}

For any element $g$ of a unital $C^*$-algebra $\CC$, we usually denote by $\sp(g)$ the spectrum
of $g$ (if $\CC$ does not have a unit we adjoin one and compute the spectrum in the canonically extended $C^*$-algebra). When precision is needed we also specify the $C^*$-algebra by writing $\sp(g\,|\,\CC)$. For
instance, if $\Sigma'$ is an invariant subset of $\Sigma$ and $g$ is a function in
$\C^\infty(\Sigma')=\CC^\infty(\Sigma')$, one has two compact subsets of $\mathbb C$:
$\sp\big(g\,|\,\CC(\Sigma')\big)$ (difficult to compute) and
$\sp\big(g\,|\,\C(\Sigma')\big)=\overline{g(\Sigma')}$. For operators $H$ in the Hilbert space $\H$, we
stick to the usual notation $\sp(H)$.

We recall a basic fact: the image of $g\in\CC$ by a unital $C^*$-morphism $\pi:\CC\to\BB$ satisfy the spectral inclusion $\,\sp\big(\pi(g)\,|\,\BB\big)\subset\sp\big(g\,|\,\CC\big)\,$, and the
spectrum is preserved if $\pi$ is injective (which is for instance the case if $\CC$ is a
$C^*$-subalgebra of $\BB$ and $\pi$ is the canonical inclusion). When $\CC,\BB$ are not unital, we can apply this remark to their minimal unitalizations.

\begin{Theorem}\label{colin}
Let $f\in\CC(\Si)$, pick $\sigma_1,\sigma_2\in\Sigma$ and set $H_{\sigma_1}:=\Op_{\sigma_1}(f)$,
$H_{\sigma_2}:=\Op_{\sigma_2}(f)$. Then
\begin{enumerate}
\item[(i)] If $\sigma_1,\sigma_2$ belong to the same orbit, then the operators
$H_{\sigma_1},H_{\sigma_2}$ are unitarily equivalent, and thus have the same spectrum (multiplicity included).

\item[(ii)] If $E_{\sigma_1}\subset E_{\sigma_2}$, then
$\sp(H_{\sigma_1})\subset\sp(H_{\si_2})\subset\sp(f\,|\,\CC(\Si))$. So, if $\sigma_1$ and $\sigma_2$ generate
the same quasi-orbit (\;\!\ie $E_{\sigma_1}=E_{\sigma_2}$)\,, then $\sp(H_{\sigma_1})=\sp(H_{\sigma_2})$.
\end{enumerate}
\end{Theorem}

\begin{proof}
(i)  Assume first that $f\in\C^\infty(\Si)$. One has $\si_2=\Theta_Z(\sigma_1)$ for some $Z\in\Xi$, which implies that
$$
f\circ\Theta_{\sigma_2}
=f\circ\Theta_{\Theta_Z(\sigma_1)}
=f\circ\Theta_{\sigma_1}\circ{\mathcal T}_Z.
$$
Therefore it is sufficient to show that $\Op(\varphi)$ and $\Op(\varphi\circ{\mathcal T}_Z)$
are unitarily equivalent if $\varphi$ belongs to the subspace $\BC^\infty(\Xi)$. But this is a well-known fact
\cite[Prop.~2.13]{Fol89}: translation by $Z$ of functions on the ``phase space" $\Xi$ leads to
unitarily equivalent Weyl quantized operators. The unitary operator $U_Z$ realizing the equivalence is the Weyl quantization $U_Z=\Op(\mathfrak e_Z)$ of the exponential
$$
\mathfrak e_Z(X):=e^{-i[\![X,Z]\!]},\ \ \ \ \ \forall X\in\Xi.
$$
Then, if we have $\,\Op_{\si_2}(f)=U_Z^*\,\Op_{\si_1}(f)\,U_Z\,$  for all $f\in\CC^\infty(\Si)$, we also get it for all $f\in\CC(\Si)$ by continuity.

\medskip
(ii) It is clear that $\sp(H_{\sigma_j})\subset\sp(f\,|\,\CC(\Si))$, because $H_{\si_j}$ is obtained from $f\in\CC(\Si)$ by applying the morphism $\Op_{\si_j}$. Actually, to have equality, it is enough that $\Op_{\si_j}$ be faithful, which happens if and only if the orbit generated by $\si$ is dense.

One can write $\Op_{\sigma_j}=\Op\circ\PP'_{\sigma_j}\circ\RR_{E_{\sigma_j}}$, with
$E_{\sigma_j}$ the quasi-orbit generated by $\sigma_j$, $j=1,2$. Since $\Op$ and
$\PP'_{\sigma_j}$ are monomorphisms, they preserve spectra; thus we only need to compare the spectrum of
$\RR_{E_{\sigma_1}}(f)$ in $\CC\left(E_{\si_1}\right)$ with the spectrum of $\RR_{E_{\sigma_2}}(f)$ in $\CC\left(E_{\si_2}\right)$. For simplicity, write $E_j$ instead of $E_{\si_j}$.

Since $E_{1}\subset E_{2}$, we define the obvious restriction mapping $\R_{21}:\C(E_{2})\rightarrow\C(E_{1})$; it is an epimorphism satisfying $\R_{E_1}=\R_{21}\circ\R_{E_2}$. By Rieffel quantization, one gets an epimorphism $\RR_{21}:\CC(E_{2})\rightarrow\CC(E_1)$ satisfying $\RR_{E_1}=\RR_{21}\circ\RR_{E_2}$. Thus $\RR_{E_{1}}(f)$ is the image of $\RR_{E_{2}}(f)$ through the morphism $\RR_{21}$ and
clearly this finishes the proof.
\end{proof}

The second part of the Theorem tells us that the natural index set for the spectra of the family $\{H_\si\}_{\si\in\Si}$ is not $\Si$, not even the orbit space $\Si/\Th$, but the smaller quotient $\Si/\!\sim$\,, where $\si\sim\si'$ means that $\si$ and $\si'$ generate the same quasi-orbit.

The next result is mainly a preparation for Theorem \ref{mailatata}.

\begin{Proposition}\label{baca}
Let $\{F\in\mathbb Q\}$ be a covering with quasi-orbits of $\,\Si$. For any $f\in\CC(\Si)$ one has
\begin{equation}\label{ighin}
\sp(f\,|\,\CC(\Si))=\overline{\bigcup_{F\in\mathbb Q}\sp\left(\RR_F(f)\,|\,\CC(F)\right)}\,.
\end{equation}
\end{Proposition}

\begin{proof}
Let us consider the morphism
$$
\RR:\CC(\Si)\rightarrow\bigoplus_{F\in\mathbb Q}\CC^F(\Si),\ \ \ \ \ \RR(f):=\left\{\RR_F(f)\right\}_{F\in\mathbb Q}
$$
(in the definition of the direct sum we require an uniform norm-bound on the family of elements, in order to have an obvious $C^*$-structure).
If we show that this morphism is injective the proof will be finished, because the spectrum of an element of a direct sum $C^*$-algebra is the closure of the union of spectra of all the components.

To show injectivity, we need to prove that
\begin{equation}\label{helf}
\bigcap_{F\in\mathbb Q}\ker\RR_F=\bigcap_{F\in\mathbb Q}\CC^F(\Si)=\{0\}.
\end{equation}
But since $\Si=\bigcup_{F\in\mathbb Q}F$ we have
\begin{equation}\label{fer}
\bigcap_{F\in\mathbb Q}\ker\R_F=\bigcap_{F\in\mathbb Q}\C^F(\Si)=\{0\},
\end{equation}
and (\ref{helf}) follows from (\ref{fer}), considering the dense common subset of smooth vectors.
\end{proof}

Of course, one can replace in Proposition \ref{baca} the space $\Si$ by any of its closed invariant subsets $\Gamma$. When $\Gamma=E$ is a quasi-orbit, one gets easily an operator version:

\begin{Corollary}\label{vaca}
Let $\{F\in\mathbb Q\}$ be a covering with quasi-orbits of a quasi-orbit $E\in\mathbf Q(\Si)$. Choose $\si\propto E$ and $\si(F)\propto F$ for every $F\in\mathbb Q$ and set $H_\si:=\Op_\si(f)$ and $H_{\si(F)}:=\Op_{\si(F)}(f)$. Then
\begin{equation}\label{then}
\sp(H_\si)=\overline{\bigcup_{F\in\mathbb Q}\sp(H_{\si(F)})}\,.
\end{equation}
\end{Corollary}

\subsection{Essential spectra}\label{etratata}

We recall the disjoint decomposition $\sp(H)=\sp_{\rm d}(H)\sqcup\spe(H)$ of the spectrum of a
(bounded) operator $H$ into its discrete and essential parts. The points $\lambda\in\sp_{\rm d}(H)$
of the discrete spectrum are, by definition, finitely degenerated eigenvalues, isolated from the rest
of $\sp(H)$. It will often be used in the sequel that the essential spectrum $\spe(H)$ of
$H\in\mathbf B(\H)$ coincides with the spectrum of the image of $H$ in the Calkin algebra
${\bf B}(\H)/{\bf K}(\H)$, where ${\bf K}(\H)$ is the (two-sided, closed) ideal of compact operators.

\begin{Proposition}\label{back}
Let $f\in\CC(\Sigma)$, choose $\sigma\in\Sigma$ of the second kind, and set $H_\sigma:=\Op_\sigma(f)$.
Then $\sp_{\rm d}(H_\sigma)=\varnothing$.
\end{Proposition}

\begin{proof}
The operator $H_\si$ belongs to $\Op\left(\BB_\si\right)$; it is enough to show that $\Op\left(\BB_\si\right)$ contains no compact operator except $0$.

We have $\B_\sigma\cap\C(\Xi)=\{0\}$ by hypothesis. Considering the smooth vectors under the action $\mathcal T$, we get
$$
\{0\}
=\big[\B_\sigma\cap\C(\Xi)\big]^\infty
=\B_\sigma^\infty\cap\C(\Xi)^\infty
=\BB_\sigma^\infty\cap\CC(\Xi)^\infty
=\big[\BB_\sigma\cap\CC(\Xi)\big]^\infty,
$$
which implies that $\{0\}=\BB_\sigma\cap\CC(\Xi)$ by density. This, together
with the injectivity of $\Op$, gives $\Op(\BB_\sigma)\cap{\bf K}(\H)=\{0\}$, which concludes the
proof.
\end{proof}

A simple picture emerges on minimal orbits:

\begin{Corollary}\label{kinder}
If the quasi-orbit $E$ is minimal, all the operators $\{H_\sigma\}_{\sigma\in E}$ have the same
spectrum, which coincides with $\sp\big(\RR_E(f)\mid\CC(E)\big)$. If $E=\Xi$, then $H_\sigma$ is a compact
operator in $\H$ and $0$ is the only point which can belong to its essential spectrum. In the opposite
case the spectrum of $H_\si$ is purely essential.
\end{Corollary}

\begin{proof}
All the points of a minimal orbit generate it, so by Theorem \ref{colin}, (ii) the operators $\{H_\si\}_{\si\in\Si}$ are equi-spectral.

If $E=\Xi$, then $H_\si\in\Op[\CC(\Xi)]=\mathbf K(\H)$. For a compact operator it is
well-known \cite[Thm.~VI.16]{RSI} that the spectrum is composed of finitely-degenerated eigenvalues
which have $0$ as the single possible accumulation point.

If $E$ is minimal but different from $\Xi$, it is of the second kind by Lemma \ref{nimalu}; then we apply Proposition \ref{back}.
\end{proof}

\begin{Remark}
{\rm When $\Si$ is compact, by a simple application of Zorn's Lemma, there always exists at least one
compact minimal quasi-orbit. So there will always be points $\sigma$ (almost periodic) for which the
operator $H_\si$ is purely essential. An extreme case is when $\Sigma$ is only composed of almost
periodic points. This happens exactly when $\Sigma$ is a disjoint union of minimal quasi-orbits and it
is equivalent to the fact that all the elements $f\in\C(\Sigma)$ are almost periodic functions (the
$\Xi$-orbit of $f$ in $\C(\Sigma)$ is relatively compact in the uniform topology). In such a situation all the operators
$H_\sigma$ have void discrete spectrum.}
\end{Remark}

We turn to a more interesting situation. From now on, we are going to denote by $f_{\Sigma'}=\R_{\Sigma'}(f)\in\C^\infty(\Sigma')$ the restriction of
the function $f\in\C^{\infty}(\Sigma)$ to the invariant subset $\Sigma'\subset\Sigma$.

\begin{Theorem}\label{mailatata}
Let $f\in\C^\infty(\Si)$ and let $\si\in\Si$ be of the first kind.
Denote by $E$ be the quasi-orbit generated by $\si$ in the
dynamical system $(\Si,\Theta,\Xi)$. For each $F\in\mathbf Q_0(E)$, choose $\si(F)\propto F$. Then
\begin{equation}\label{mailatatu}
{\rm sp}_{{\rm ess}}\left[H_\si\right]={\sp}\left[\,f_{E^{{\rm n}}}\mid\CC(E^{{\rm n}})\right]=
\overline{\bigcup_{F\in\mathbf Q_0(E)}{\rm sp}\left[f_F\mid\CC(F)\right]}=
\overline{\bigcup_{F\in\mathbf Q_0(E)}{\rm sp}\left[H_{\si(F)}\right]}.
\end{equation}
\end{Theorem}

\begin{proof}
The last equality follows from the fact that
$$
\CC(F)\ni f_F\mapsto H_{\si(F)}=\Op\left[\PP'_\si(f_F)\right]\in\mathbf B(\H)
$$
is a monomorphism of $C^*$-algebras. The second equality is a consequence of Proposition \ref{baca}. Just replace $\Si$ with $E^{{\rm n}}$ and take $\mathbb Q=\mathbf Q_0(E)=\mathbf Q(E^{{\rm n}})$.

We are going to justify the first equality, using the $\Xi$-morphisms introduced above.
We know that ${\rm sp}_{{\rm ess}}\left[H_\si\right]$ equals the usual spectrum of the image of $H_\si$ into
$$
\Op\left[\BB_\si\right]/\mathbf K(\H)\cong\BB_\si/\CC(\Xi)\cong\CC(E)/(\PP'_\si)^{-1}[\CC(\Xi)].
$$
By functoriality, this last quotient is the quantization of $\C(E)/(\P'_\si)^{-1}[\C(\Xi)]$. By Proposition \ref{detai}, one has $(\P'_\si)^{-1}[\C(\Xi)]=\C^{E^{{\rm n}}}(E)$, so it follows that
$$
\C(E)/(\P'_\si)^{-1}[\C(\Xi)]\cong\C(E^{{\rm n}}),
$$
implying (once again by functoriality) that
$$
\CC(E)/(\PP'_\si)^{-1}[\CC(\Xi)]\cong\CC\left(E^{{\rm n}}\right).
$$
It will follow that $\CC\left(E^{{\rm n}}\right)$ is isomorphic to $\Op\left[\BB_\si\right]/\mathbf K(\H)$. By inspection, it is easy to see that $f_{E^n}$ is send by this isomorphism into the image of $H_\si$ in this quotient and so the first equality in (\ref{mailatatu}) is proven.
\end{proof}

\begin{Remark}\label{porcu}
{\rm The proof of the Theorem shows that if $\QQ_1(E)\subset\QQ_0(E)$ such that $E^{{\rm n}}=\bigcup_{F\in \QQ_1(E)}F$,
then we also have
$$
{\rm sp}_{{\rm ess}}\left[H_\si\right]=\overline{\bigcup_{F\in\mathbf Q_1(E)}{{\rm sp}}\left[f_F\mid\CC(F)\right]}=
\overline{\bigcup_{F\in\mathbf Q_1(E)}{\rm sp}\left[H_{\si(F)}\right]}.
$$
The operators $H_{\si(F)}$ could be called {\it asymptotic Hamiltonians} for $H_\si$. If one of these asymptotic Hamiltonians is null (equivalent to $f_F=0$ for some $F\in\QQ_0(E)$), then the point $0$ belongs to the essential spectrum of $H_\si$.}
\end{Remark}

\begin{Remark}
{\rm We stated Theorem \ref{mailatata} only for $f\in\C^\infty(\Si)$ in order to have simpler notations. It clearly extends (with the same proof) to all $f\in\CC(\Si)$; just replace the restrictions $f_F$ and $f_{E^{{\rm n}}}$ by $\RR_F(f)$ and $\RR_{E^{{\rm n}}}(f)$, respectively.}
\end{Remark}

\subsection{Some examples}\label{seci}

The results of the previous sections are general; they apply to any $\Xi=\mathbb R^{2n}$-dynamical system. To get concrete examples, it would be nice to understand at least partially the quasi-orbit structure of $\Si$. This can be achieved in many cases and it seems to be pointless to draw a large list; the reader can try his own particular cases. In \cite{ABG,Da,GI1,GI2,GI3,HM,LS,Ma1,MPR2} one encounters many examples of configuration space anisotropy (connected to the behavior of the symbol $f$ in the variable $x\in\X$) which can be adapted to full phase-space anisotropy; of course the results will be different. We are going to indicate briefly only a couple of interesting instances, stressing the advantages inherent to the present setting. The examples will be centered around Theorem \ref{mailatata} and will involve only smooth functions, for simplicity.

Although very particular, the following situation covers zero order pseudodifferential operators with phase-space anisotropic symbols: We assume that $\Si$ is a compactification of $\Xi$, \ie it is a compact space containing $\Xi$ as a dense open subset (maybe after an identification). The action $\Th$ of $\Xi$ on $\Si$ is a continuous extension of the action $\mathcal T$ of $\Xi$ on itself by translations. The dense orbit $\Xi$ is the set of generic points of $\Si$, while the closed complement $\Si_\infty:=\Si\setminus\Xi$ consists of non-generic points. We are going to regard the Abelian $C^*$-algebra $\mathcal A=\C(\Si)$ directly as the $\Xi$-algebra of all continuous functions on $\Xi$ which can be extended continuously on $\Si$. Then, plainly, the elements $f$ of $\C^\infty(\Si)$ are those which are smooth on $\Xi$ and for which all the derivatives have this extension property. One is interested in $H:=H_{\si=0}=\Op(f)$ (all the other $H_X$\,, for $X\in\Xi\subset\Si$, are unitarily equivalent to this one). The operators $H_\si$ with $\si\in\Si_\infty$ are only used to express the essential spectrum of $H$. As a consequence of Theorem \ref{mailatata} one can write
$$
{\rm sp}_{{\rm ess}}\left[H\right]=
\overline{\bigcup_{F}{\rm sp}\left[H_{\si(F)}\right]},
$$
where for each quasi-orbit $F\subset\Si_\infty$ a generating point $\si(F)$ was chosen.

\begin{Example}\label{shat}
{\rm A simple non-trivial particular case is formed of {\it vanishing oscillation functions} $f\in{\rm VO}(\Xi)$. These are complex continuous functions on $\Xi$ such that for any compact subset $K$ of $\Xi$ one has
$$
\left[{\rm osc}_K(f)\right](X):=\underset{Y\in K}{\sup}|f(X+Y)-f(X)|\underset{X\rightarrow\infty}{\longrightarrow}0.
$$
It is shown easily that it is a $\Xi$-algebra and that its spectrum $\Si$ can be identified to a compactification of $\,\Xi$ such that all the elements of $\Si_1:=\Si\setminus\Xi$ are fixed points under the extension $\Th$ of the action by translations $\mathcal T$. This is the largest example for which "the quasi-orbits at infinity" are reduced to points. Then it follows easily that the "asymptotic Hamiltonians" $\{H_\si\}_{\si\in\Si_1}$ are just the constant operators $\Op(c)=c\,{\rm id}_{L^2(\X)}$ constructed with all the values $c$ taken by the function $f$ at infinity. Consequently the essential spectrum of $H$ coincides with {\it the asymptotic range} of $f$:
$$
{\rm sp}_{{\rm ess}}\left[H\right]=R_{{\rm asy}}(f):=\bigcap_{K\in\kappa(\Xi)}\overline{f(\Xi\setminus K)},
$$
where $\kappa(\Xi)$ is the family of all the compact neighborhoods of the origin in $\Xi$.}
\end{Example}

\begin{Example}\label{shut}
{\rm  It is known that the $C^*$-algebra ${\rm AP}(\Xi)$ of all continuous almost periodic functions on $\Xi$ is $\Xi$-simple. Equivalently, its Gelfand spectrum $b\Xi$ ({\it the Bohr group associated to} $\Xi$) is a minimal dynamical system. This would lead immediately to the absence of discrete spectrum for large classes of almost periodic pseudodifferential operators. We consider more interesting to mix this class with ${\rm VO}(\Xi)$. Let us denote by $\left<{\rm VO}(\Xi)\cdot{\rm AP}(\Xi)\right>$ the smallest $C^*$-algebra containing both ${\rm VO}(\Xi)$ and ${\rm AP}(\Xi)$. It is a $\Xi$-algebra with spectrum $\Si=\Xi\sqcup(\Si_1\times b\Xi)$. The non-generic quasi-orbits have all the form $\{\si_1\}\times b\Xi$ for some $\si_1\in\Si_1$. For a smooth element $f$ of $\left<{\rm VO}(\Xi)\cdot{\rm AP}(\Xi)\right>$, the essential spectrum of the operator $H$ can be written in terms of purely almost periodic operators. A simple very explicit case is $f=gh$, with $g\in{\rm VO}(\Xi)^\infty$ and $h\in{\rm AP}(\Xi)^\infty$. One has
$$
\sp_{{\rm ess}}[\Op(gh)]=\overline{\bigcup_{\si_1\in\Si_1}\sp\left[\Op(g(\si_1)h)\right]}=
\overline{\bigcup_{\si_1\in\Si_1}g(\si_1)}\,\sp[\Op(h)]=R_{{\rm asy}}(g)\sp[\Op(h)].
$$
In the same way one shows that, under the same assumptions on $g$ and $h$, we get
$$
\sp_{{\rm ess}}[\Op(g+h)]=R_{{\rm asy}}(g)+\sp[\Op(h)].
$$
Extensions of this result to many classes of minimal functions, generalizing $AP(\Xi)$, are available by the approach of \cite{Ma3}.
}
\end{Example}

Other type of compactifications, which are quite different from Examples \ref{shat} and \ref{shut}, are suggested by the decomposition $\Xi=\X\times\X^*$. If $\Omega$ is a compactification of $\Xi$ and $\Omega^*$ a compactification of $\X^*$, then $\Si:=\Omega\times\Omega^*$ will be a compactification of $\Xi$ and $\C(\Si)\cong\C(\Omega)\otimes\C(\Omega^*)$. It is natural to consider actions $\Th=\th\otimes\th^*$, where $\th$ is an action of $\X$ on $\Omega$ extending the translations in $\X$ and $\th^*$ is an action of $\X^*$ on $\Omega^*$ extending the translations in $\X^*$. We are in a position to apply the results above. We leave to the reader the task to write down quasi-orbits for this situation and to make statements about essential spectra. We are only going to outline a situation contrasting to Example \ref{shat}.

\begin{Example}\label{shit}
{\rm We consider the $\Xi$-algebra $\mathcal A={\rm VO}(\X)\otimes{\rm VO}(\X^*)$, where ${\rm VO}(\X)$ and ${\rm VO}(\X^*)$ are defined in an obvious way. Its Gelfand spectrum can be written as $\Si=(\X\sqcup\Omega_1)\times(\X^*\sqcup\Omega_1^*)$, where the points of $\Omega_1$ are fixed by the action $\th$ extending the translations on $\X$ and analogously for $\Omega^*_1$. Aside the big quasi-orbit $\Si=\Omega\times\Omega^*$, one still has other types of quasi-orbits defined by points $\omega\in\Omega_1$ and $\omega^*\in\Omega_1^*$:
\begin{enumerate}
\item
$\{\omega\}\times\Omega^*$ generated by $(\omega,\xi)$ for any $\xi\in\X^*$,
\item
$\Omega\times\{\omega^*\}$ generated by $(x,\omega)$ for any $x\in\X$,
\item
$\{(\omega,\omega^*)\}$.
\end{enumerate}
The first two types will suffice, because of Remark \ref{porcu}. For any smooth element $f$ in ${\rm VO}(\X)\otimes{\rm VO}(\X^*)$, setting $H:=\Op_{\si=0}(f)$, we get
$$
\sp_{{\rm ess}}(H)=\overline{\bigcup_{\omega\in\Omega_1}\sp\left(f|_{\{\omega\}\times\Omega^*}\mid\CC(\{\omega\}\times\Omega^*)\right)}
\,\bigcup\,\overline{\bigcup_{\omega^*\in\Omega^*_1}
\sp\left(f|_{\Omega\times\{\omega^*\}}\mid\CC(\Omega\times\{\omega^*\})\right)}.
$$
It is easy to see that $\CC\left(\Omega\times\{\omega^*\}\right)$ is isomorphic to the Abelian $C^*$-algebra $\C(\Omega)$ while $\CC\left(\{\omega\}\times\Omega^*\right)$ is isomorphic to the Abelian $C^*$-algebra $\C(\Omega^*)$. This leads straightforwardly to
$$
\sp_{{\rm ess}}(H)=f(\Omega_1\times\Omega^*)\cup f(\Omega\times\Omega^*_1).
$$
We note that {\it the quantum quadrant}, contained in Chapter 11 of \cite{Rie1}, is a $C^*$-subalgebra of the quantization of $\mathcal A={\rm VO}(\mathbb R)\otimes{\rm VO}(\mathbb R^*)$, thus it is covered by our treatment.

One might also want to work out the case $\mathcal A=\left<{\rm VO}(\X)\cdot{\rm AP}(\X)\right>\otimes\left<{\rm VO}(\X^*)\cdot{\rm AP}(\X^*)\right>$. On the other hand, decompositions of $\Xi$ in direct sums different from $\X\times\X^*$ can also lead to interesting situations.}
\end{Example}

\begin{Remark}\label{remark}
{\rm The attentive reader might have observed that the previous Example is built on the short exact sequence
$$
0\rightarrow\C(\X)\otimes\C(\X^*)\rightarrow{\rm VO}(\X)\otimes{\rm VO}(\X^*)\rightarrow \left[\C(\Omega_1)\otimes{\rm VO}(\X^*)\right]\oplus\left[{\rm VO}(\X)\otimes\C(\Omega_1^*)\right]\rightarrow 0.
$$
To help, we notice that $\C(\Xi)\cong\C(\X)\otimes\C(\X^*)$ and that ${\rm VO}(\X)/\C(\X)\cong\C(\Omega_1)$ and ${\rm VO}(\X^*)/\C(\X^*)\cong\C(\Omega^*_1)$. Such short exact sequences can be written for all the possible tensor products, but the extra fact that the dynamics in $\Omega_1,\Omega_1^*$ are trivial helped to get explicit quasi-orbits and thus explicit contributions to the essential spectrum. However, the cornerstone was the possibility to turn the exact sequence of Abelian $C^*$-algebras into an exact sequence of non-commutative $C^*$-algebras, these ones being those concerned by the spectral analysis of the pseudodifferential operators. Both Rieffel's functor and the crossed product are exact functors, but the crossed product cannot cover most of the phase-space types of anisotropy of the symbol $f$.
}
\end{Remark}

\begin{Example}\label{shet}
{\rm We come back to example \ref{shat} and remark that ${\rm VO}(\Xi)$ contains the $C^*$-algebra $\C_{{\rm rad}}(\Xi)$ of all the continuous functions admitting radial limits at infinity, i.e those which can be extended to the radial compactification $\Si_{{\rm rad}}:=\Xi\sqcup S^{2n-1}$ (obvious topology; the points at infinity will be fixed points).
Thus our results apply easily. In \cite{Rie1} the $C^*$-algebra $\CC(\mathbb R^2\sqcup S^1)$ is called {\it the quantum euclidean closed disk}.

Following Chapter 12 in \cite{Rie1} (see also references therein), one can inflate this example by gluing together several discs $\mathbb R^2\sqcup S^1$ along the circle $S^1$ into a dynamical system which is no longer topologically transitive. The circle $S^1$ will be composed of fixed points and the interior of each disc will be an orbit. Theorem \ref{mailatata} applies easily and with an explicit output to all the Hamiltonians given by the points of the dynamical system. {\it The quantum sphere} is obtained with a pair of discs. To get {\it a quantum version of the group} $SU(2)$ one uses a family of discs parametrised by the one-dimensional torus.
}
\end{Example}

\begin{Example}
{\rm A simple example which is not topologically transitive and for which the action $\Th$ is not the extension of some translations is the following: We let $\Xi=\X\times\X^*=\mathbb R\times\mathbb R$ act on itself by $\Th_{(x,\xi)}(y,\eta):=\left(e^xy,e^\xi\eta\right)$; the abelian $C^*$-algebra is $\A:=\C(\Xi)$. Since the action is not given by translations, the quantized version $\mathfrak A$ will no longer be elementary (i.e. isomorphic to the ideal of all the compact operators in an infinite-dimensional, separable Hilbert space).

There are obviously nine quasi-orbits: four closed quarter-planes, four coordinate semi-axes (all containing the origin) and the origin, which is a fixed point. The points of the open quarter-planes are of the first kind and all the others are of the second kind. The generic and the non-generic points in each orbit are evident. Denoting as usual $H_\si:=\Op(f\circ\Th_\si)$ for some smooth element $f$ of $\mathfrak A$ and for all the points $\si=(y,\eta)\in\Si=\Xi$, one gets by Theorem \ref{mailatata}
$$
\sp_{{\rm ess}}\left[H_{(1,1)}\right]=\sp\left[H_{(1,0)}\right]\cup\sp\left[H_{(0,1)}\right],\ \ \ \ \
\sp_{{\rm ess}}\left[H_{(1,-1)}\right]=\sp\left[H_{(1,0)}\right]\cup\sp\left[H_{(0,-1)}\right],
$$
$$
\sp_{{\rm ess}}\left[H_{(-1,1)}\right]=\sp\left[H_{(-1,0)}\right]\cup\sp\left[H_{(0,1)}\right],\ \ \ \ \
\sp_{{\rm ess}}\left[H_{(-1,-1)}\right]=\sp\left[H_{(-1,0)}\right]\cup\sp\left[H_{(0,-1)}\right].
$$
On the other hand, since the points belonging to the semi-axes are of the second type, the operators $H_{(\pm1,0)}$ and $H_{(0,\pm1)}$ are purely essential. It is a simple exercise to work out the structure of $\CC(\mathbb R_+\times\{0\})$ (it is abelian), to show that $\left[H_{(1,0)}u\right](y)=f(e^y,0)\,u(y)$ for all $u\in L^2(\mathbb R)$, so $\sp\left[H_{(1,0)}\right]=\overline{f(\mathbb R_+\times\{0\})}$. We have analogous results for the other semi-axes, so the essential spectra of all the possible operators $H_\si$ are known explicitly. Clearly $H_{(0,0)}=f(0,0)\,{\rm id}_{L^2(\mathbb R)}$.

This Example goes under the name {\it the algebraists' real quantum plane} in \cite{Rie1}, Chapter 12. Obvious higher-dimensional instances are available.}
\end{Example}

\subsection{Random operators}\label{secintre}

The framework studied above appears rather often in the context of random families of operators, where
some extra structure is present. The essence of the theory of random operators is to study families
$\{H_\sigma\}_{\sigma\in\Sigma}$ of self-adjoint operators indexed by a set $\Si$ on which a
probability measure $\mu$ is given, being mainly interested in properties that hold with probability
one. In most cases the probability measure is invariant under the ergodic action of a group and
the family of operators has an equivariance property.

We place this general idea in the framework introduced so far. Let us assume
that the dynamical system $(\Sigma,\Theta,\Xi)$ is compact and metrisable, and that $\Sigma$ is endowed
with a $\Theta$-invariant and ergodic probability measure $\mu$ defined on a $\si$-algebra including all the open sets.
Recall that ergodicity means that the
$\Theta$-invariant subsets of $\Sigma$ must have measure zero or one. The metrisability condition on
$\Sigma$ implies that $\CC(\Si)$ is separable. The notion of {\em hull} of a physical system leads to
such a setting (\cite{BHZ}).

\begin{Proposition}\label{pastur}
As before, let $f\in\CC(\Si)$ and set $H_\si=\Op_\si(f)$ for any $\si\in\Si$.
\begin{enumerate}
\item
There exists a closed set $S\subset\mathbb R$ such that $\sp(H_\sigma)=S\,$ for $\mu$-almost every
$\sigma\in\Sigma$.
\item
One has $\sp_{\rm d}(H_\sigma)=\varnothing\,$ for $\mu$-almost every $\sigma\in\Sigma$.
\end{enumerate}
\end{Proposition}

\begin{proof}
1. We recall that the topological support $\supp(\mu)$ of the probability measure $\mu$ is the
smallest closed set $M\subset\Sigma$ such that $\mu(M)=1$. Now we know from Lemma 3.1
of \cite{BHZ} (based on Birkhoff's Ergodic Theorem) that $\supp(\mu)$ is a quasi-orbit and the set
$$
\Si_0:=\{\si\in\Si\mid E_\si=\supp(\mu)\}
$$
is measurable and $\mu(\Si_0)=1$.
So the claim follows by Proposition \ref{colin} (ii).

\medskip
2. Due to Proposition \ref{back}, it is sufficient to show that there exists a measurable set
$\Sigma_1\subset\Sigma$, with $\mu(\Sigma_1)=1$, such that $\B_\sigma$ is of the second kind for each
$\sigma\in\Sigma_1$.

Let us denote by $E$ the closed set $\supp(\mu)$, which is the quasi-orbit generated by the
points $\sigma\in\Sigma_0$, with $\mu(\Sigma_0)=1$.
Then $(E,\Theta,\mu,\Xi)$ is once again a compact, metrisable, ergodic dynamical system. By
Birkhoff's Ergodic Theorem, there exists a measurable set $\Sigma_1\subset\Sigma_0\subset E$, with
$\mu(\Sigma_1)=1$, such that for each $\sigma\in\Sigma_1$
$$
\int_E\d\mu(\si')\,g(\sigma')
=\lim_{R\to\infty}\frac1{|B_R|}\int_{B_R}\d X\,g\big[\Theta_X(\sigma)\big].
$$
Here $B_R:=\{X\in\Xi\mid|X|\le R\}$ and $g$ is any positive element of $\C(E)$. If
$g\circ\Theta_\sigma\in\C(\Xi)$, we obtain immediately
$$
\lim_{R\to\infty}\frac1{|B_R|}\int_{B_R}\d X\,g\big[\Theta_X(\sigma)\big]=\lim_{R\to\infty}\frac1{|B_R|}\int_{B_R}g\circ\Th_\si=0.
$$
Thus
$$
\int_E\d\mu(\si')\,g(\sigma')=0,
$$
which implies $g=0$. Doing this for the positive and negative parts of (respectively) the real and the imaginary part
of an arbitrary element $f\in\C(E)$, it will follow that $f=0$ as soon as $f\circ\Theta_\sigma\in\B_\sigma\cap\C(\Xi)$. So
$\B_\sigma\cap\C(\Xi)=\{0\}$, and thus $\B_\sigma$ is of the second kind for each
$\sigma\in\Sigma_1$.
\end{proof}

Such results any many others (almost sure constancy of the spectral types of the family $\{H_\si\}_{\si\in\Si}$) are classical in the theory of ergodic random families of operators. They rely on the equivariance condition expressed in Theorem \ref{colin}, (i) and can be proven in a more abstract framework, using only measurability assumptions \cite{CL,PF}. We included Proposition \ref{pastur} here because the proof fits nicely in our setting and because the statement can be made somewhat more precise as usually:
The proof supplies an explicit example (not unique, of course) of a set of full measure for which the corresponding family of spectra is constant: $\Si_0=\supp(\mu)^{{\rm g}}$ is the set of generic points of the quasi-orbit $\supp(\mu)$.
On the other hand, using the decomposition $\,\Si=\Si_I\sqcup\Si_{II}\,$ in subsets of points of the first, respectively second kind, the preceding proof shows that $\Si_I$ is $\mu$-negligible. This and Proposition \ref{back} are the reasons for having purely essential operators with probability one.

\subsection{The semiclassical limit of spectra}\label{serata}

In Quantum Mechanics one also encounters {\it the Planck constant} $\hb$ which has been conventionally taken equal to $1$ until now. We let it vary in the interval $(0,1]$ and study the continuity of spectra of the emerging operators.

When $\hb$ is taken into account, in the formula for the usual Weyl quantization one has to replace (\ref{op}) by
\begin{equation}\label{oph}
\left[\mathfrak{Op}^\hb(h)u\right](x):=(2\pi\hb)^{-n}\int_\X\int_{\X^*}dx\,d\xi\,e^{\frac{i}{\hb}(x-y)\cdot
\xi}\,h\left(\frac{x+y}{2},\xi\right)u(y),
\end{equation}
which also requires replacing the composition law (\ref{rodact}) on $\C^\infty(\Si)$ by
\begin{equation}\label{rodacti}
f\,\#^\hb\,g:=(\pi\hb)^{-2n}\int_\Xi\int_\Xi dYdZ\,e^{\frac{2i}{\hb}[\![Y,Z]\!]}\,\Theta_Y(f)\,\Theta_Z(g)=
\end{equation}
$$
=\pi^{-2n}\int_\Xi\int_\Xi dYdZ\,e^{2i[\![Y,Z]\!]}\,\Theta_{\sqrt\hb Y}(f)\,\Theta_{\sqrt h Z}(g).
$$

The entire formalism works exactly as for the case $\hb=1$ (cf. \cite{Rie1}, where somewhat different notations and conventions are used) and for each $\hb$ one gets a quantized $C^*$-algebra $\CC^\hb(\Si)$ (with composition law $\#^\hb$ and norm $\parallel\cdot\parallel_{\CC^\hb(\Si)}$) having the same properties and allowing the same constructions as $\CC(\Si)\equiv\CC^{\hb=1}(\Si)$. Even the spectral results above have their obvious $\hb$-counterparts; this will be used below.

In addition, in \cite{Rie1} it is shown that the family $\left\{\CC^\hb(\Si)\right\}_{\hb\in[0,1]}$ can be organized in a continuous field of $C^*$-algebras; it actually provides a strict deformation quantization of a natural Poisson algebra constructed on $\C^\infty(\Si)$. It is obvious from (\ref{rodacti}) that  the $C^*$-algebra $\CC^\hb(\Si)$ is obtained by applying the general procedure to the classical data $\big(\Si,\Th^{\hb},\Xi,[\![\cdot,\cdot]\!]\big)$, where $\Th^{\hb}_X:=\Th_{\sqrt\hb X}$ for any $(\hb,X)$. The $C^*$-algebra $\CC^{\hb=0}(\Si)$ is simply taken to be $\C(\Si)$. We remark that $\C^\infty(\Si)$ is a dense $^*$-subalgebra of any $\CC^\hb(\Si)$.

Exactly as before and using (\ref{oph}), for every $\hb\in(0,1]$ and every $\si\in\Si$, we can construct the representations
$$
\Op^\hb_\si:\CC^\hb(\Si)\rightarrow\mathbf B(\H),\ \ \ \ \ \Op^\hb_\si:=\Op^\hb\circ\PP^\hb_\si,
$$
which can be used to supply families of $\hb$-quantum Hamiltonians.

Sending to \cite{Rie1} for details, we outline now only the facts that will be used in the proof of Theorem \ref{miclasica}. The classical data $\left([0,1]\times\Si,\Th',\Xi,[\![\cdot,\cdot]\!]\right)$ can also be considered, where $\Th'_X(\hb,\si):=\left(\hb,\Th^\hb_X(\si)\right)$ for every $X,\hb,\si$. This gives raise by quantization to the $C^*$-algebra $\CC([0,1]\times\Si)$. Now we take into account the $\Xi$-epimorphisms
$$
\mathcal N^\hb:\C([0,1]\times\Si)\rightarrow\C(\Si),\ \ \ \ \ \left[\mathcal N^\hb(\mathfrak f)\right](\si):=\mathfrak f(\hb,\si).
$$
Since they intertwines the actions $\Th'$ and $\Th^\hb$, they are send by the Rieffel functor into epimorphisms
$$
\mathfrak N^\hb:\CC([0,1]\times\Si)\rightarrow\CC^\hb(\Si),\ \ \ \ \ \mathfrak N^\hb|_{\C^\infty([0,1]\times\Si)}=\mathcal N^\hb|_{\C^\infty([0,1]\times\Si)}.
$$
A basic fact, contained in the definition of a continuous field and proven in \cite{Rie1}, is that the mapping
$$
[0,1]\ni\hb\mapsto\parallel\mathfrak N^\hb(\mathfrak g)\parallel_{\CC^\hb(\Si)}
$$
is continuous for any $\mathfrak g\in\CC([0,1]\times\Si)$.

Clearly, $\Si$ can be replaced by any closed invariant subset $\Gamma$ in all the considerations above. In the proof of Theorem \ref{miclasica} we are going to take $\Gamma=E\in\mathbf Q(\Si)$.

\medskip
After all these preparations, let us introduced the concept of continuity for families of sets that will be useful in Theorem \ref{miclasica}.

\begin{Definition}
 Let $I$ be a compact interval and suppose given a family $\left\{S^\hb\right\}_{\hb\in I}$ of closed subsets of
 $\,\mathbb{R}$.
\begin{enumerate}
\item The family $\left\{S^\hb\right\}_{\hb\in I}$ is called {\rm
outer continuous} if for any $\hb_0\in I$ and any compact subset $K$ of $\,\mathbb R$ such that
$K\cap S^{\hb_0}=\varnothing$, there exists a neighborhood $V$ of $\hb_0$ with
$K\cap S^{\hb}=\varnothing$, $\forall \hb\in V$.
\item The family $\left\{S^\hb\right\}_{\hb\in
I}$ is called {\rm inner continuous} if  for any $\hb_0\in I$ and any open subset $A$ of $\,\mathbb R$ such
that $A\cap S^{\hb_0}\ne\varnothing$, there exists a neighborhood $W\subset I$ of
$\hb_0$ with $A\cap S^{\hb}\ne\varnothing$, $\forall \hb\in W$.
\item If the family is both inner and outer continuous, we say simply that it is {\rm continuous}.
\item Sometimes, to express continuity at a point $\hb_0\in I$, we write suggestively $S^\hb\rightarrow S^{\hb_0}\,$ for $\hb\rightarrow\hb_0$.
\end{enumerate}
\end{Definition}

In the proof of the next result we are going to use the functional calculus for self-adjoint operators, so we shall ask the function $f$ to be real.

\begin{Theorem}\label{miclasica}
For any $f\in\C^\infty(\Si)_\mathbb R$, $\si\in\Si$ and $\hb\in(0,1]$ we set $H^\hb_\si:=\Op^\hb_\si(f)$ and $S^\hb_\si:=\sp\left(H^\hb_\si\right)$. For $\hb=0$ we set $S^0_\si:=\overline{f(E_\si)}$. Then the family of compact sets $\left\{S^\hb_\si\right\}_{\hb\in[0,1]}$ is inner and outer continuous. In particular one has
$$
\sp\left(H^\hb_\si\right)\rightarrow \overline{f\left(E_\si\right)}\ \ \ \ {\rm when}\ \hb\rightarrow 0.
$$
\end{Theorem}

\begin{proof}
(i) First we recall that the conclusion of the Theorem follows if it is proven that the mapping
\begin{equation}\label{cuci}
[0,1]\ni\hb\mapsto\left\Vert\left(H^\hb_\si-\zeta\right)^{-1}\right\Vert_{\mathbf B(\H)}
\end{equation}
is continuous for any $\zeta\notin\mathbb R$. This is Proposition 2.5 in \cite{MP3} (see also the references therein). The proof is straightforward, it also works for unbounded self-adjoint operators and we shall not repeat it here.

\medskip
(ii) Let us denote by $g^{(-1)_\hb}$ the inverse of $g$ with respect to the composition law $\#^\hb$. Then (\ref{cuci}) follows if we show that
$$
[0,1]\ni\hb\mapsto\left\Vert\left(\RR^\hb_{E_\si}(f)-\zeta\right)^{(-1)_\hb}\right\Vert_{\CC^\hb(E_\si)}
$$
is continuous for any $\zeta\notin\mathbb R$, since the standard Weyl representation $\Op^\hb$ is faithful. Since $f$ has been chosen to be a smooth vector, one has $\RR^\hb_\si(f)=f_{E_\si}$. In addition, there is no loss of generality to assume that $E_\si=\Si$. Thus we are reduced to show for any $f\in\C^\infty(\Si)$ that
$$
[0,1]\ni\hb\mapsto\left\Vert\left(f-\zeta\right)^{(-1)_\hb}\right\Vert_{\CC^\hb(\Si)}
$$
is continuous for any $\zeta\notin\mathbb R$.

\medskip
(iii) We define
$$
\mathfrak f(\hb,\si):=f(\si),\ \ \ \ \ \forall\,(\hb,\si)\in[0,1]\times\Si.
$$
Obviously $\mathfrak f$ belongs to $\C^\infty([0,1]\times\Si)\subset\CC([0,1]\times\Si)$ and it is a self-adjoint element. Let us introduce $\mathfrak r_\zeta:=(\mathfrak f-\zeta)^{[-1]}$, where $\,^{[-1]}$ indicates inversion in $\CC([0,1]\times\Si)$; notice that we have $\mathfrak N^\hb(\mathfrak r_\zeta)=(f-\zeta)^{(-1)_\hb}$.

By the continuous field property, the mapping
$$
[0,1]\ni\hb\mapsto\left\Vert\mathfrak N^\hb(\mathfrak r_\zeta)\right\Vert_{\CC^\hb(\Si)}=\left\Vert(f-\zeta)^{(-1)_\hb}\right\Vert_{\CC^\hb(\Si)}
$$
is continuous and this finishes the proof.
\end{proof}

Combining previous results one also gets the semiclassical limits of essential spectra, in the setting of Theorem \ref{miclasica}.

\begin{Corollary}\label{niclasica}
\begin{enumerate}
\item
If $\si$ is of the first kind, then $\,\sp_{{\rm ess}}\left(H^\hb_\si\right)\rightarrow\overline{f(E_\si^{{\rm n}})}\,$ when $\hb\rightarrow 0$.
\item
If $\si$ is of the second kind, then $\,\sp_{{\rm ess}}\left(H^\hb_\si\right)\rightarrow\overline{f(E_\si)}\,$ when $\hb\rightarrow 0$.
\end{enumerate}
\end{Corollary}

\begin{proof}
Assertion $2\,$ follows from Theorem \ref{miclasica} and the obvious extension of Proposition \ref{back} to arbitrary $\hb\in(0,1]$, saying that $\sp_{{\rm d}}\left(H_\si^\hb\right)=\emptyset$ if $\si$ is of the second kind.

The point $1\,$ is a consequence of Theorem \ref{miclasica} and a simple adaptation of Theorem \ref{mailatata}.
\end{proof}

\bigskip
\bigskip
{\bf Acknowledgements:} The idea of writing this article originated in discussions with Vladimir Georgescu and Eduardo Friedman. We are grateful to Rafael Tiedra de Aldecoa for his interest in the project.

The author is partially supported by {\it N\'ucleo Cientifico ICM P07-027-F "Mathematical
Theory of Quantum and Classical Magnetic Systems"} and by Chilean Science Foundation {\it Fondecyt} under the Grant 1085162.


\bigskip
\bigskip
\bigskip
{\bf Address}

\medskip
Departamento de Matem\'aticas, Universidad de Chile, 

Las Palmeras 3425, Casilla 653, Santiago, Chile

\emph{E-mail:} Marius.Mantoiu@imar.ro

\end{document}